\documentclass[10pt]{article}
\usepackage{amssymb}
\begin{document}
\title{The 2-color Rado Number of $x_1+x_2+\cdots +x_{m-1}=ax_m,\  \textrm{II}$}
\author{Dan Saracino\\Colgate University}
\date{}
\maketitle
\begin{abstract} In the first installment of this series, we proved that, for every integer $a\geq 3$ and every $m\geq 2a^2-a+2$, the 2-color Rado number of $$x_1+x_2+\cdots +x_{m-1}=ax_m$$ is $\lceil\frac{m-1}{a}\lceil\frac{m-1}{a}\rceil\rceil.$  
 Here we obtain the best possible improvement of  the bound on $m.$  We prove that if $3|a$ then the 2-color Rado number is $\lceil\frac{m-1}{a}\lceil\frac{m-1}{a}\rceil\rceil$  when $m\geq 2a+1$ but not when $m=2a,$ and that if $3\nmid a$ then the 2-color Rado number is $\lceil\frac{m-1}{a}\lceil\frac{m-1}{a}\rceil\rceil$ when $m\geq 2a+2$ but not when $m=2a+1.$  We also determine the 2-color Rado number for all $a\geq 3$ and $m\geq \frac{a}{2}+1.$ 
\end{abstract}
\vspace{.25in}

\noindent \textbf{1. Introduction}

\vspace{.25in}

A special case of the work of Richard Rado [\textbf{5}] is that for every integer $m\geq 3$ and  all positive integers $a_1,\ldots,a_m$ there exists a smallest positive integer $n$ with the following property:  for every coloring of the elements of the set $[n]=\{1,\ldots,n\}$ with two colors, there  exists a solution of the equation $$a_1x_1+a_2x_2+\cdots +a_{m-1}x_{m-1}=a_mx_m$$ using elements of $[n]$ that are all colored the same. (Such a solution is called \emph{monochromatic}.) The integer $n$ is called the \emph{2-color Rado number} of the equation.

In 1982, Beutelspacher and Brestovansky [\textbf{1}] proved that for every $m\geq 3$, the 2-color Rado number of $$x_1+x_2+\cdots +x_{m-1}=x_m$$ is $m^2-m-1.$  In 2008 Guo and Sun [\textbf{2}] generalized this result by proving that, for all positive integers $a_1,\ldots,a_{m-1},$ the 2-color Rado number of the equation $$a_1x_1+a_2x_2+\cdots+a_{m-1}x_{m-1}=x_m$$  is $aw^2+w-a,$ where $a=\textrm{min}\{a_1,\dots, a_{m-1}\}$ and $w=a_1+\cdots +a_{m-1}.$
In the same year, Schaal and Vestal [\textbf{7}] dealt with the equation $$x_1+x_2+\cdots +x_{m-1}=2x_m.$$ They proved, in particular, that for every $m\geq 6,$ the 2-color Rado number is $\lceil\frac{m-1}{2}\lceil\frac{m-1}{2}\rceil\rceil.$  Building on the work of Schaal and Vestal, we proved in $[\textbf{6}]$ that for every $a\geq3$ and $m\geq 2a^2-a+2,$ the 2-color Rado number of the equation $x_1+\cdots +x_{m-1}=ax_m$ is $\lceil\frac{m-1}{a}\lceil\frac{m-1}{a}\rceil\rceil.$   Our main purposes here are to obtain the best possible improvement of the bound on $m,$ and to determine the Rado number in most cases where $m$ falls below the improved bound.

We begin by using a sharpening  of the arguments in [\textbf{6}] to prove (in Section 3) the following result.

\vspace{.15in}

\noindent\textbf{Theorem 1}.  For every integer $a\geq 3$ and every $m\geq a^2-a+1,$ the 2-color Rado number of the equation$$x_1+x_2+\cdots +x_{m-1}=ax_m$$ is $\lceil\frac{m-1}{a}\lceil\frac{m-1}{a}\rceil\rceil.$

\vspace{.15in}

\noindent\textbf{Notation.}  We will denote $\lceil\frac{m-1}{a}\lceil\frac{m-1}{a}\rceil\rceil$ by $C(m,a),$ and we will denote the equation indicated in the statement of Theorem 1 by $L(m,a).$  We will denote the 2-color Rado number of $L(m,a)$ by$R_2(m,a).$

\vspace{.15in}

In order to present the rest of our results efficiently, we next prove (in Section 4) the following.

\vspace{.15in}

\noindent\textbf{Theorem 2.}   Suppose $a+1\leq m\leq 2a+1.$ Then $R_2(m,a)=1$ iff  $m=a+1.$  If $a+2\leq m\leq 2a+1,$ then $R_2(m,a)\in \{3,4,5\},$ and we have:

\begin{itemize}

\item[] $R_2(m,a)=3$  iff $m\leq \frac{3a}{2}+1$ and $a\equiv m-1$ (mod 2).

\item[] $R_2(m,a)=4$ iff either:

\begin{itemize}
\item[(i)] $m\leq \frac{3a}{2}+1$ and $a\not \equiv m-1$ (mod 2), or

\item[(ii)] $m> \frac{3a}{2}+1$ and $a\equiv m-1$  (mod 3).

\end{itemize}

\item[] $R_2(m,a)=5$ iff $m> \frac{3a}{2}+1$ and $a\not \equiv m-1$  (mod 3).

\end{itemize}

\vspace{.15in}

Theorem 2 will be useful to us in Section 5, where we obtain our final  lowering of the bound on $m,$  which is as follows.

\vspace{.15in}

\noindent\textbf{Theorem 3.}  Suppose $a\geq 3.$ If $3|a$ then $R_2(m,a)=C(m,a)$ when $m\geq 2a+1$ but $R_2(2a,a)=5$ and $C(2a,a)=4.$ If $3\nmid a$ then $R_2(m,a)=C(m,a)$  when $m\geq 2a+2$ but $R_2(2a+1,a)=5$ and $C(2a+1,a)=4.$

\vspace{.15in}

By the results of [\textbf{7}], Theorem 3 also holds when $a=2.$

Finally, in Section 6, we prove Theorems 4 and 5, which determine all values of $R_2(m,a)$  when $\frac{a}{2}+1\leq m\leq a.$

\vspace{.15in}

\noindent\textbf{Theorem 4.}  If  $\frac{2a}{3}+1\leq m\leq a,$ then:

\begin{itemize} \item[] for $a=3$ we have $R_2(a,a)=9,$ and \item[] for  $a\geq 4$ we have

\begin{itemize} \item[] $R_2(m,a)=3$ if $a\equiv m-1$ (mod 2) and \item[] $R_2(m,a)=4$ if $a \not \equiv m-1$  (mod 2).  

\end{itemize}
\end{itemize}

\vspace{.15in}

\noindent\textbf{Theorem 5.}  If $\frac{a}{2}+1\leq m < \frac{2a}{3}+1$ (so $a\geq 4)$  then:

\begin{itemize}
\item[] for $a\equiv m-1$  (mod 3) we have $R_2(m,a)=4,$ and

\item[] for $a\not \equiv m-1$ (mod 3) we have $R_2(m,a)=5$ \emph{except} that
\begin{itemize}
\item[]  $R_2(3,4)=10$ and $R_2(4,5)=9,$ and

\item[] $R_2(m,a)=6$ if $10\leq a\leq 14$ and $m=a-4.$

\end{itemize}

\end{itemize}

\vspace{.15in}

\noindent\textbf{Conventions and definitions.} In working with a fixed 2-coloring of $[n]$, we will use the colors red and blue, and we will denote by $R$ and $B,$ respectively, the sets of elements colored red and blue.  We will call a 2-coloring of $[n]$ \emph{bad} if it yields no monochromatic solution of $L(m,a).$

\vspace{.25in}

\noindent\textbf{2. Preliminary lemmas}

\vspace{.25in}

The results of $[\textbf{6}]$ relied on the fact that if $m\geq 2a^2-a+2$ then $2m-2\leq C(m,a),$ and therefore numbers in $[2m-2]$ can be used in producing solutions of $L(m,a)$ in $[C(m,a)].$  The improvement presented in Theorem 1 rests on showing that we can obtain the same results using  $[m-1]$ instead of $[2m-2]$, and that $[m-1]\subseteq [C(m,a)]$ if $m\geq a^2-a+2.$  (The case $m=a^2-a+1$ will be handled separately, in Proposition 1 below.)

\vspace{.15in}

 \noindent\textbf{Lemma 1.}   Suppose $a\geq 3$ and $m\geq a^2-a+2.$  Then  $m-1\leq C(m,a).$  

\vspace{.15in}

\noindent\emph{Proof.}  If $m\geq a^2+1,$ then $\frac{(m-1)^2}{a^2}\geq m-1 $ and the result follows. If $a^2-a+2\leq m\leq a^2$  we can write $m=a^2-a+b,$ where $2\leq b\leq a.$  We have $\frac{m-1}{a}=a-1+\frac{b-1}{a},$
and therefore $\lceil \frac{m-1}{a} \rceil=a$ and $C(m,a)=a^2-a+b-1=m-1.
$ $\Box$
 
 \vspace{.15in}

It is shown in Proposition 1 of $[\textbf{6}]$ that, for $m\geq 3,$  $C(m,a)$ is a lower bound for $R_2(m.a).$  So to prove Theorem 1 we must show that, for $m\geq a^2-a+1,$ $C(m,a)$ is also an upper bound, i.e., every 2-coloring of $[C(m,a)]$ yields a monochromatic solution of $L(m,a).$

To proceed, it will be convenient to recall the  compact notation used in $[\textbf{6}]$ to indicate solutions of $L(m,a).$

\vspace{.15in}

\noindent\textbf{Notation.}  If $n_1,\ldots,n_k$ are nonnegative integers whose sum is $m$, and $d_1,\ldots,d_k$ are elements of $[C(m,a)]$ such that we obtain a true equation from $L(m,a)$ by substituting $d_1$ for the variables $x_1,\ldots, x_{n_1}$, $d_2$ for the next $n_2$ variables, and so on, then we denote this true equation by $$[n_1\rightarrow d_1;\ n_2\rightarrow d_2;\ \cdots;\  n_k\rightarrow d_k].$$

\vspace{.15in}

For example, the true instance 
$$
 a+a+\cdots +a=a(m-1)
$$
 of $L(m,a)$ will be denoted by $$[m-1\rightarrow a; \ 1\rightarrow m-1].$$

 \vspace{.15in}

\noindent \textbf{Proposition 1.}  If $m=a^2-a+1,$ then every 2-coloring of $[C(m,a)]$ yields a monochromatic solution of $L(m,a).$

\vspace{.15in}

\noindent \emph{Proof.}  Note that if $m=a^2-a+1,$ then $C(m,a)=(a-1)^2.$

Suppose we have a bad 2-coloring of $C(m,a),$ and suppose, without loss of generality, that $1\in R.$ Then the solution $[a^2-a\rightarrow 1;\ 1\rightarrow a-1]$ shows that $a-1\in B,$ and multiplying the assigned values  in this solution by $a-1$ shows that $(a-1)^2\in R.$ But the solution $[(a-1)^2\rightarrow 1;\ a\rightarrow (a-1)^2]$ shows that $(a-1)^2\in B,$ a contradiction.  $\Box$ 

\vspace{.15in}

By Proposition 1, we can assume, in completing the proof of Theorem 1, that $
m\geq a^2-a+2,$ and therefore Lemma 1 applies.

Some of our arguments in Section 3 will require $a\geq 4.$  When $a=3,$ Theorem 1 asserts that $R_2(m,3)=C(m,3)$ for $m\geq 7,$ and this is proved in Section 6 of $[\textbf{6}]$.  Accordingly, we need only consider $a\geq 4$ in what follows.

\vspace{.15in}

\noindent\textbf{Conventions.}  In the remainder of Section 2, and in Section 3, we assume that $a \geq 4$ and $m\geq a^2-a+2.$  We suppose that we have a bad 2-coloring of $[C(m,a)],$  and we seek a contradiction. We assume without loss of generality that $a-2\in R.$

\vspace{.15in}

As in $[\textbf{6}]$, we proceed by considering two cases, depending on the coloring of the element $a-1.$  If $a-1\in B,$ then we can obtain our contradiction by using the same argument as in $[\textbf{6}]$, since that  argument uses only elements in $[m-1].$
(See $[\textbf{6}]$, Section 3.)  Accordingly, we adopt another convention.

\vspace{.15in}

\noindent\textbf{Convention.}  We assume in the remainder of Section 2, and in Section  3,  that $a-1\in R.$

\vspace{.15in}

\noindent\textbf{Lemma 2.}  The elements $1$ and $a$ are in $R$.

\vspace{.15in}

The proof is as in Lemmas 4 and 5 of $[\textbf{6}],$ which use only numbers in $[m-1].$

\vspace{.15in}

 \noindent\textbf{Lemma 3.}  The numbers $m-a,\ldots, m-1$ are all in $ B.$
 
 \vspace{.15in}
 
 \noindent\emph{Proof.}  We want  to show that $m-a+j\in B$ for $0\leq j\leq a-1.$ Since  $1, \ a-1,$ and $a$ are all in $R$ and we are assuming that there are no monochromatic solutions of $L(m,a)$ in $[C(m,a)],$ we need only consider the solution $$[m-2a+2j+1\rightarrow a;\ a-1-j\rightarrow a-1;\ a-1-j\rightarrow 1;\ 1\rightarrow m-a+j.]\ \ \Box$$  
  
  \vspace{.15in}
  
  \noindent\textbf{Lemma 4.}  The numbers $1,2,\ldots, a$ are all in $R$. 
  
  \vspace{.15in}
  \noindent\emph{Proof}.  For $0\leq j\leq a-1,$ consider the solution
  
   $$[m-a+j\rightarrow j+1;\ a-j\rightarrow m-a+j]$$ and use the result of Lemma 3.
  $\Box$
  
  \vspace{.15in}
  
The next result generalizes Lemma 9 from $[\textbf{6}].$

\vspace{.15in}

\noindent\textbf{Lemma 5.}  If $d$ is an integer such that $a|d$ and $m-1\leq d\leq a(m-1)$, then $\frac{d}{a}\in B.$

\vspace{.15in}

\noindent\emph{Proof.}  Write $d=(m-1)j+k,$ with $1\leq j\leq a-1$ and $0\leq k\leq m-1.$  Then the solution
$$ \left[m-1-k\rightarrow j;\ k\rightarrow j+1;\ 1\rightarrow \frac{d}{a}\right] $$ shows that $\frac{d}{a}\in B.$  $\Box$

\vspace{.25in}

\noindent\textbf{3.  The proof of Theorem 1}

\vspace{.25in}

In this section we will use the results of Section 2,  together with algebraic expressions for $C(m,a),$ to produce a red solution of $L(m,a)$ in $[C(m,a)]$.  This will contradict our standing assumption that our 2-coloring of $[C(m,a)]$ is bad, and conclude the proof of Theorem 1.

The following Lemma is Lemma 10 from $[\textbf{6}].$

\vspace{.15in}

\noindent\textbf{Lemma 6.}  Let $m=ua^2+va+c,$ with $u$ as large as possible and $0\leq v,c\leq a-1.$
\begin{itemize}
\item[(i)]  If $c=1$ then $C(m,a)=\frac{(m-1)^2}{a^2}.$
\item[(ii)]  If $c=0$ then $C(m,a)=\frac{m^2-m+va}{a^2}.$
\item[(iii)] If $2\leq c\leq a-1$ then $C(m,a)=\frac{m^2+(a-c-1)m+c-ac-vac+va+ta^2}{a^2},$ \\where $t=\left\lceil\frac{(c-1)(v+1)}{a}\right\rceil.$
\end{itemize}

\vspace{.15in}

When $c=1,$ the argument in $[\textbf{6}]$ produces a red solution of $L(m,a)$ by using only elements of $C(m,a)$ that can be shown to be in $R$ by using elements of $[m-1]$.  So the same argument yields a red solution here.  We turn to the remaining cases.

\vspace{.15in}

\noindent  \textbf{The Case $c=0$}

\vspace{.15in}

In this case we have $\frac{m}{a}\in B$ by Lemma 5. We choose an $s$ such that $s\in R, \ s+1\in B,$ and $s+1\leq \frac{m}{a}.$  Using the expression for $C(m,a)$ in Lemma 6,  we obtain $$C(m,a)=\left(\frac{m-a}{a}\right)\frac{m}{a}+\frac{(a-1)m+va}{a^2}.$$ We let $$\alpha=\frac{m-a}{a}(s+1)+\frac{(a-1)m+va}{a^2}\leq \left(\frac{m-a}{a}\right)\frac{m}{a}+\frac{(a-1)m+va}{a^2},$$ so $\alpha\leq C(m,a).$  As in $[\textbf{6}],$ we see that $\alpha\in R.$

 We now obtain a red solution of $L(m,a)$  by assigning the value $\alpha$ to $x_{m-2},x_{m-1}$ and $x_m$ and the value $s$ to $(a-2)(\frac{m-a}{a})$ other variables, and showing that we can assign values in $R$ to the remaining $\frac{2m}{a}+a-5$ variables to complete the solution.  In fact we will show that we can accomplish this by using only values in the set $[a]$. These values are all in $R$ by Lemma 4.

The values assigned to the remaining variables must add up to
$$\frac{a-2}{a}(m-a)+\frac{a-2}{a^2}((a-1)m+va).$$

If we can show that using only the value $a$ yields a sum that is at least this large, and using only the value $1$ yields a sum that is at most this large, then there is a unique solution that uses  values  in one of the sets $\{j,\ j+1\}$, where $j\in [a-1].$

Since $v\leq a-1,$ we can achieve our first objective by showing that 
$$a\left(\frac{2m}{a}+a-5\right)\geq \frac{a-2}{a}(m-a)+\frac{a-2}{a^2}((a-1)m+(a-1)a),$$  which simplifies to $$a^2-5a+1-\frac{2}{a}\geq m\left(\frac{2-5a}{a^2}\right).$$   Since the right-hand side is negative, this is clearly true when $a\geq 5.$  When $a=4$  it is true since $m\geq 14$ because $m\geq a^2-a+2.$

Since $v\geq 0,$ we can achieve our second objective by showing that
$$ \left(\frac{2m}{a}+a-5\right)\leq \frac{a-2}{a}(m-a)+\frac{a-2}{a^2}((a-1)m).$$  But this simplifies to $2a^3-7a^2\leq m(2a^2-7a+2),$ which is true for all $a\geq 4$ and $m\geq a.$ (It is \emph{not} true when $a=3$ and $m> 9,$ and this is why we dealt with the case $a=3$ separately at the outset.)
 
\vspace{.15in}

\noindent \textbf{The Case $2\leq c\leq a-1$}

\vspace{.15in}

In this case we have $\frac{m+a-c}{a}\in B$ by Lemma 5. We choose an $s$ such that $s\in R,\ s+1\in B,$ and $s+1\leq \frac{m+a-c}{a}.$  Using the expression for $C(m,a)$ in Lemma 6, we obtain $$
C(m,a)=\left(\frac{m-c}{a}\right)\left(\frac{m+a-c}{a}\right)+\frac{(c-1)(m+a-c) +a\gamma}{a^2},$$ where $\gamma=ta-(c-1)(v+1)$, with  $t$ as in Lemma 6.  Note that since

$$0\leq t-\frac{(c-1)(v+1)}{a}\leq 1$$ by definition of $t$, we have $0\leq \gamma\leq a.$

  We consider the element
$$
\beta=\left(\frac{m-c}{a}\right)(s+1)+\frac{(c-1)(m+a-c)+a\gamma}{a^2}\leq C(m,a).$$
To see that $\beta\in R,$ we consider  the solution $$\left[m-c\rightarrow s+1;\ c-2\rightarrow \frac{m+a-c}{a};\ 1\rightarrow \frac{m+a-c}{a}+\gamma;\ 1\rightarrow \beta\right].$$  Note that $\frac{m+a-c}{a}+\gamma\in B$ by Lemma 5, since  $$\frac{m+a-c}{a}+\gamma\leq \frac{m+a-c+a^2}{a}\leq \frac{m+a-2+a^2}{a}$$ and it is easy to verify that $m+a-2+a^2\leq a(m-1)$ when $a\geq 4$ and $m\geq a^2-a+2.$

To obtain our  red solution of $L(m,a)$, we assign the value $\beta$ to $x_m, x_{m-1}$ and $x_{m-2},$ and the value $s$ to $(a-2)(\frac{m-c}{a})$ other variables, and show that we can assign values in $R$ to the remaining $\frac{2(m-c)}{a}+c-3$ variables to complete the solution.  We again use values in the set $[a].$

The values assigned to the remaining $\frac{2(m-c)}{a}+c-3$ variables must add up to 
\begin{equation} \frac{a-2}{a}(m-c)+\frac{a-2}{a^2}((c-1)(m+a-c)+a\gamma). \end{equation}   If we can show that using only the value $a$ (respectively, $1$) yields a sum that is at least (respectively, at most) this large, then, as before, there must be a solution that uses values in one of the sets $\{j,\ j+1\}$, where $j\in [a-1].$ 
  
   Using the fact that $\gamma\leq a,$  we can achieve our first objective by showing that
$$a\left(\frac{2(m-c)}{a}+c-3\right) \geq \frac{a-2}{a}(m-c)+\frac{a-2}{a^2}((c-1)(m+a-c)+a^2), $$  which simplifies to $$c^2(a-2)+c(a^3-2a^2-a+2)+(-4a^3+3a^2-2a)\geq m(-a^2-3a+2+c(a-2)).$$  If we regard $a$ as a constant and denote the quantity on the left-hand side of this inequality by $f(c)$, then the derivative $$f'(c)=2c(a-2)+(a^3-2a^2-a+2)$$ is easily seen to be positive for $c\geq 0$ and $a\geq 4,$ so the minimum value of  $f(c)$ for $2\leq c \leq a-1$ occurs at $c=2.$ Since $$m(-a^2-3a+2+c(a-2))\leq m(-a^2-3a+2+(a-1)(a-2))=m(4-6a),$$ we only need to verify that $f(2)\geq m(4-6a),$ i.e., $$2a^3+ a^2+4\leq m(6a-4), $$ and this is  true for $a\geq 4$ and $m\geq a^2-a+2.$

To achieve our second objective, it will suffice, by using expression (1) and the fact that $0\leq \gamma$, to show that
$$\left(\frac{2(m-c)}{a}+c-3\right) \leq \frac{a-2}{a}(m-c)+\frac{a-2}{a^2}((c-1)(m+a-c)).$$  This inequality simplifies to
$$c^2(a-2)+c(a^2-3a+2)-2a^2-2a\leq m(a^2-5a+2+c(a-2)).$$  Denoting the quantity on the left-hand side by $g(c)$, we have $$g'(c)=2c(a-2)+(a^2-3a+2),$$ so $g'(c)>0$ for $c\geq 0$ and $a\geq 4.$  Therefore the maximum value of $g(c)$ for $2\leq c \leq a-1$ occurs at $c=a-1.$ Since 
$$m(a^2-5a+2+c(a-2))\geq m(a^2-5a+2+2(a-2))=m(a^2-3a-2),$$
we need only verify that $g(a-1)\leq m(a^2-3a-2),$ i.e., that $$2a^3-10a^2+8a-4\leq m(a^2-3a-2).$$  This is easily verified for $a\geq 4$ and $m\geq a^2-a+2.$ (But it fails when $a=3$ and $m> 8,$ again indicating why we dealt separately with the case $a=3$.)  $\Box$

\vspace{.25in}

\noindent\textbf{4. The proof of Theorem 2}

\vspace{.25in}

The following lemma will be useful in proving Theorems 2,\ 4, and 5.

\vspace{.15in}

\noindent\textbf{Lemma 7.}  Suppose $\frac{a}{2}+1\leq m \leq 2a+1.$ Then:

\begin{itemize}

\item[(1)]  there exists a solution of $L(m,a)$ using only values in $\{1,2\},$
\item[(2)] there exists a solution of $L(m,a)$ using only values in $\{1,3\}$ iff \\$a\equiv m-1$  (mod 2),

\item[(3)] there exists a solution of $L(m,a)$ using only values in $\{2,3\}$ iff \\$\frac{2a}{3}+1\leq m\leq \frac{3a}{2}+1,$ and

\item[(4)] there exists a solution of $L(m,a)$ using only values in $\{1,4\}$ iff \\$a\equiv m-1$  (mod 3).

\end{itemize}

\vspace{.15in}

\noindent\emph{Proof.}  There exists a solution using values in $\{1,2\}$ iff either
$m-1\leq a\leq 2(m-1)$ or $m-1\leq 2a\leq 2(m-1).$ If $\frac{a}{2}+1\leq m\leq a+1$ then the first alternative holds, and if $a+1\leq m\leq 2a+1$ then the second holds. This proves statement (1).  There exists a solution of $L(m,a)$ using values in $\{2,3\}$ iff either $2(m-1)\leq2a\leq 3(m-1)$ or $2(m-1)\leq 3a\leq 3(m-1),$ i.e., iff $\frac{2a}{3}\leq m-1\leq \frac{3a}{2}.$  This proves statement (3).

As we assign values in $\{1,3\}$ to $x_1,
\ldots,x_{m-1},$ the total values achieved by the left side of $L(m,a)$ are exactly those integers that have the same parity as $m-1$ and are between $m-1$ and $3(m-1),$ inclusive. So there exists a solution using values in $\{1,3\}$ iff $a\equiv m-1$  (mod 3) and either $m-1\leq a\leq 3(m-1)$ or $m-1\leq 3a\leq 3(m-1).$  As in the proof of (1), one of these pairs of inequalities must hold, and this proves (2).  The proof of (4) is similar.  $\Box$

\vspace{.15in}

\noindent\emph{The proof of Theorem 2.}  Assume that $a+1\leq m\leq 2a+1.$

It is clear that $R_2(m,a)=1$ iff we can obtain a solution of $L(m,a)$ by assigning all the variables the same value, and this is so iff $m-1=a.$  Note that $R_2(m,a)$ can never be 2, for if we color 1 and 2 differently then we can only obtain a monochromatic solution of $L(m,a)$ in $[2]$ by coloring all the variables the same, but then $R_2(m,a)=1.$

For the remainder of the proof we assume $a+2\leq m\leq 2a+1,$ so that $R_2(m,a)\geq 3,$ and there is no solution of $L(m,a)$ that assigns all the variables the same color.

We next establish the conditions under which $R_2(m,a)=3.$  Suppose first that we have a bad 2-coloring of $[3]$, with, say, $1\in R.$   By statement (1) of Lemma 7, there is a solution of $L(m,a)$ using values in $\{1,2\}$, and both values must be used n the solution. So $2\in B.$ It then follows from statements (2) and (3) of Lemma 7 that if $a\equiv m-1$  (mod 2) and $m\leq \frac{3a}{2}+1$ then we must color 3 both blue and red in order to avoid a monochromatic solution of $L(m,a)$ in $[3],$ so $R_2(m,a)=3.$ If  $a\not \equiv m-1$ (mod 2) (or, respectively, if $m> \frac{3a}{2}+1$) then we can color 3 red (or, respectively, blue) and obtain a bad 2-coloring of $[3],$ so $R_2(m,a)> 3.$

By what we have just shown, either of conditions (i) or (ii) in the statement of Theorem 2 implies that $R_2(m,a)\geq 4.$  To prove that each implies $R_2(m,a)=4,$ we suppose we have a bad 2-coloring of $[4]$, with $1\in R,$ and we seek a contradiction, assuming that (i) or (ii) holds.  By using statement (1) of Lemma 7, and then doubling all the values in its proof, we get $2\in B$ and $4\in R.$

If (i) holds then by statement (3) of Lemma 7 we have $3\in R.$  To get a contradiction, we obtain a red solution of $L(m,a).$ We start by assigning the value 4 to $x_m$ and the value 1 to each of the other variables. We must show that we can increase the total value of the left side by $4a-(m-1)$ by increasing some of the values on the left side by 2 or 3. So we need to write $4a-(m-1)$ as a sum of 2's and 3's, using at most $m-1$ terms. Our bounds on $m$ imply that $2\leq 4a-(m-1)\leq 3(m-1),$ so this is possible (using at most two 2's).

If (ii) holds then by statement (4) of Lemma 7 we can obtain a red solution of $L(m,a)$ using values in $\{1,4\}.$

We have shown that each of (i) and (ii) implies $R_2(m,a)=4.$  Conversely, if $R_2(m,a)=4$ then since $R_2(m,a)\neq 3,$ either (i) holds or $m> \frac{3a}{2}+1.$  In the latter case we must have $a\equiv (m-1)$ (mod 3), for otherwise by statements (3) and (4) of Lemma 7, the coloring $R=\{1,4\}, B=\{2,3\}$ is bad, contradicting $R_2(m,a)=4.$

Finally, suppose $m>\frac{3a}{2}+1$ and $a\not \equiv m-1$ (mod 3).  Then $R_2(m,a)\geq 5.$ To prove equality, suppose for a contradiction that we have a bad 2-coloring of $[5]$, with $1\in R.$ As above, we see that $2\in B$ and $4\in R,$ and as in our proof that condition (i) implies $R_2(m,a)=4$ we get a contradiction if $3\in R.$  So suppose $3\in B.$

We claim that $5\in R.$ To see this we construct a solution of $L(m,a)$ in which we assign the value 5 to $x_m$ and values in $\{2,3,5\}$ to all the other variables. If we start by assigning the value 2 to all the other variables, then we must increase the value of the left side by $5a-2(m-1)$ by increasing some of the 2's by 1 or 3 each.  Note that $5a-2(m-1)\geq 0$ since $m\leq 2a+1,$ and, since $a\leq m-2,$ we have $5a-2(m-1)\leq 5(m-2)-2(m-1)=3m-8.$  Any nonnegative integer less than or equal to $3m-5$ can be written in the form $3q+r,$ with $0\leq q\leq m-2$ and $0\leq r\leq 2,$ with $r\leq 1$ if $q=m-2$. So we can achieve the desired solution (using at most two 1's), and $5\in R.$

We now obtain a red solution of $L(m,a)$ (and therefore a contradiction) by assigning the value 4 to $x_m$ and values in $\{1,4,5\}$ to all the other variables. If we start by assigning the value 1 to each of $x_1,\ldots, x_{m-1}$, then to finish we must write $4a-(m-1)$ as a sum of 3's and 4's, using at most $m-1$ terms.  Our bounds on $m$ yield $4a-(m-1)\geq 4a-2a\geq 6$ and $4a-(m-1)\leq 4(m-1).$  Therefore it is easy to show that the desired expression for $4a-(m-1)$ exists.

We have shown that if $m> \frac{3a}{2}+1$ and $a\not \equiv m-1$ (mod 3) then $R_2(m,a)=5.$ Conversely, if $R_2(m,a)=5$ then $R_2(m,a)$ is neither 3 nor 4, so by what we have already shown, we must have $m> \frac{3a}{2}+1$ and $a\not \equiv m-1$ (mod 3).  $\Box$
 \vspace{.25in}
 
 \noindent\textbf{5.  The proof of Theorem 3}
 
 \vspace{.25in}
 
 We have $C(2a+1,a)=4,$ and, since $a\geq 3, $ $C(2a,a)=4$ as well.  If $3|a,$ then by Theorem 2 we have $R_2(2a+1,a)=4$ and $R_2(2a,a)=5, $ while if  $3\nmid a$ then  $R_2(2a+1,a)=5.$ Therefore, to prove Theorem 3 it will suffice to prove the following.
 
 \vspace{.15in}
 
 \noindent\textbf{Proposition 2.}  For $a\geq 3$ and $m\geq 2a+2,$ $R_2(m,a)=C(m,a).$
 
 \vspace{.15in}
 
 By Theorem 2 of [\textbf{6}], we know that Proposition 2 holds when $a=3$.  By Theorem 1, we also know that it holds when $m>a^2-a.$  So we adopt the following conventions for the remainder of this section.
 
 \vspace{.15in}
 
 \noindent\textbf{Conventions.}  We have $a\geq 4,$ and $2a+2\leq m\leq a^2-a.$  We write $m=av+c,$ with $2\leq v\leq a-1$ and $0\leq c \leq a-1,$ so that when $v=2$ we have $c\geq 2$ and when $v=a-1$ we have $c=0.$  We suppose that we have a bad 2-coloring of $[C(m,a)]$, with $1\in R,$ and we seek a contradiction.
 
 \vspace{.15in}
 
 We consider three cases.
 
 \vspace{.15in}
 
 \noindent\textbf{Case 1:}  $1,2\in R$
 
 \vspace{.15in}
 
 \noindent\textbf{Lemma 8.}  When $1,2\in R,$ we have $C(m,a)\in R.$
 
 \vspace{.15in}
 
 \noindent\emph{Proof.}  Since $m\geq 2a+2,$ we have $\frac{m-1}{a}\geq2+\frac{1}{a}$ and $\lceil\frac{m-1}{a}\rceil\geq 3.$  So if we let $n=\lceil\frac{m-1}{a}\rceil$ then $n+1\leq C(m,a)$ and, for $k\in \{n,\ n+1\},$
 $$m-1\leq ak\leq 2(m-1).$$   So we may assign either of the values $n,\ n+1$ to $x_m$ and obtain a solution of $L(m,a)$ by assigning a value of 1 or 2 to each of $x_1,\ldots, x_{m-1}.$  Since $1,2\in R,$ we have $n,\ n+1\in B.$
 
 To show that $C(m,a)\in R,$ it therefore suffices to show that 
 $$n(m-1)\leq aC(m,a)\leq (n+1)(m-1).$$  The first inequality holds because $n=\lceil\frac{m-1}{a}\rceil.$  The second inequality asserts that $C(m,a)\leq n(\frac{m-1}{a})+\frac{m-1}{a},$ and this is true because $C(m,a)$ exceeds $n(\frac{m-1}{a})$ by less than 1.  $\Box$
 
 \vspace{.15in}
 
 We now obtain a contradiction by showing that, for some positive integer $j\leq a,$ we obtain a red solution of $L(m,a)$ by assigning the value $C(m,a)$ to $x_m$ and to $a-j$ of the variables $x_1,\ldots,x_{m-1},$ and assigning the value 1 or 2 to each of the remaining $a(v-1)+j+c-1$ variables.  To show this, we must show that for some positive $j\leq a,$
 
 \begin{equation} a(v-1)+j+c-1\leq jC(m,a)\leq 2(a(v-1)+j+c-1).
 \end{equation}
 
 \vspace{.1in}
 
 \noindent\emph{Subcase 1:} $c=0$
 
 \vspace{.1in}
 
 In this case, Lemma 6 yields $C(m,a)=v^2,$ so we must show that
 $$a(v-1)+j-1\leq jv^2\leq 2(a(v-1)+j-1)$$ for some positive $j\leq a.$  The first inequality clearly holds when $j=a.$  We now choose $j$ to be the smallest positive integer such that the first inequality holds. We claim that, for this $j,$ the second inequality holds.
 
 If $j=1$ the second inequality says that $v^2\leq 2a(v-1),$ and since $a\geq v+1$ it suffices to show that $v^2\leq 2(v^2-1).$  But this is clearly true, since $v\geq 2.$  If $j> 1$ then by the minimality of $j$ we have $a(v-1)+j-2> (j-1)v^2,$ so $jv^2< v^2 +a(v-1) +j-2$ and we need to show that
 $$v^2+a(v-1)+j-2\leq 2(a(v-1)+j-1).$$  This inequality reduces to $v^2\leq a(v-1)+j,$ and since $a\geq v+1$ it suffices to show that $v^2\leq v^2-1+j,$ which is clearly true. 
 
 \vspace{.1in}
 
 \noindent\emph{Subcase 2:} $c=1$
 
 \vspace{.1in}
 
 The argument for this case is nearly identical to that for $c=0.$  We omit the details.
 
 \vspace{.1in}
 
 \noindent\emph{Subcase 3: } $2\leq c\leq a-1$
 
 \vspace{.1in}
 
 In this case, Lemma 6 yields $C(m,a)=v^2+v+t,$ where $t=\lceil\frac{(c-1)(v+1)}{a}\rceil.$ So when $j=a$ the first inequality in statement (2) says that $av+c-1\leq a(v^2+v+t)$, which is clear.  We choose the smallest positive $j$ such that the first inequality of statement (2) holds.
 
 If $j=1,$ the second inequality in (2) says that $v^2+v+t\leq 2(a(v-1)+c).$  Since $c\neq 0,$ our conditions on $m$ imply that $a\geq v+2,$ so since $t\leq v+1$ it will suffice to show that $v^2+2v+1\leq 2(v^2+v-2+c).$ This reduces to $5\leq v^2+2c,$  which is clearly true.
 
 If $j> 1$ then by the minimality of $j$ we have $$a(v-1)+j+c-2> (j-1)(v^2+v+t),$$
 so $j(v^2+v+t)< a(v-1)+j+c-2+v^2+v+t,$ and to verify the second inequality in (2) we want to show that $$a(v-1)+j+c-2 +v^2+v+t\leq 2(a(v-1)+j+c-1),$$
 i.e., $v^2+v+t\leq a(v-1) +j+c.$  Since $c\neq 0$ implies $a\geq v+2,$ it suffices to show that $v^2+v+t\leq v^2+v-2+j+c,$ i.e., $t\leq j+c-2.$  But $t\leq c-1,$ so this is clear.
 
 \vspace{.15in}
 
 \noindent\textbf{Case 2:} $1\in R, \ 2\in B,\ 3\in B$
 
 \vspace{.15in}
 
 \noindent\emph{Subcase 1:} $c-1< \frac{a}{3}$
 
 \vspace{.1in}
 
 In this subcase we will produce a red solution of $L(m,a)$ by using values that are at most $3v$.  Note that if $c=0$ or 1 then $C(m,a)=v^2$ and $v\geq 3,$ so it is clear that $3v\leq C(m,a)$.  If $c\geq 2$ then $C(m,a)=v^2+v+t$ and it is again clear that $3v\leq C(m,a).$
 
 By assigning a value of 2 or 3 to each of the variables $x_1,\ldots, x_{m-1},$ we can achieve for the left side of $L(m,a)$ any total value between $2(av+c-1)$ and $3(av+c-1),$ inclusive.  By the assumption of the current subcase, this implies that when $c> 0$ we can achieve a solution of $L(m,a)$ using 2's and 3's on the left side and any of $2v+1,\ldots, 3v$ on the right.  When $c=0$ we can use any of $\{2v,\ldots,3v-1\}$ on the right.  So $2v+1,\ldots, 3v\in R$ when $c> 0$ and $2v,\ldots,3v-1 \in R$ when $c=0.$
 
 When $c=0,$ $$[m-a-1\rightarrow 1; \ a-1\rightarrow 2v; \ 1\rightarrow 2v+1; \ 1\rightarrow 3v-1]$$  is a red solution of $L(m,a).$
 
To obtain a red solution when $c> 0$, start by assigning  the value $3v$ to $x_m$ and the value 1 to each of the other variables. We must then increase the total value of the left side by $2av-(c-1),$ which is easily seen to be at least $3(2v)$ for any $a\geq 4,$ since $v\geq 2$ and $c\leq a-1.$  Write $$2av-(c-1)=q(2v)+r,$$ with $3\leq q< a$ and $0\leq r\leq 2v-1.$ If we increase the values of each of $x_1,\ldots, x_q$ to $2v+1,$ then we must still increase the total value of the left side by $r$.  Since $r\leq 2v-1$ and $q\geq 3,$ we can accomplish this by again increasing the values of some of $x_1,\ldots,x_q$ without increasing any value $2v+1$ by more than $v-1$ (even if $v-1=1$).
 
 \vspace{.1in}
 
 \noindent\emph{Subcase 2:} $c-1\geq \frac{a}{3}$
 
 \vspace{.1in}
 
 In this subcase we will use values no larger than $3v+1.$  Note that $3v+1\leq C(m,a)$ since $c\geq 2,$ so $C(m,a)=v^2+v+t$  and $t\geq 1.$
 
 By the assumption of the current subcase, we see as in the preceding subcase that we now have $2v+2,\ldots,3v+1\in R.$
 
 To obtain a red solution of $L(m,a)$ we start by assigning the value $3v+1$ to $x_m$ and the value 1 to each of the other variables. We must then increase the total value of the left side by $a(2v+1)-(c-1),$ which is easily shown to be least $4(2v+1)$ for any $a\geq 5.$   Assuming for the moment that $a\geq 5,$ write $$a(2v+1)-(c-1)=q(2v+1)+r,$$ with $4\leq q<  a$ and $0\leq r\leq 2v.$
 Increase the values of $x_1,\ldots,x_q$ to $2v+2$. Since $r\leq 2v$ and $q\geq 4,$ we can then increase the total value of the left side by $r$ by increasing some of the values $2v+2$ by no more than $v-1$ each.
 
 If $a=4,$ then by our bounds on $m$ we have $10\leq m\leq 12.$ In the current subcase we also have $c-1\geq \frac{4}{3},$ so $c\geq 3.$ Thus $m=11,$ and $[10\rightarrow 2;\ 1\rightarrow 5]$ and $[6\rightarrow 2;\ 4\rightarrow 3;\ 1\rightarrow 6]$ are solutions of $L(m,a).$ Since $2,3\in B,$ we have $5,6\in R,$ so $[8\rightarrow 1;\ 2\rightarrow 6;\ 1\rightarrow 5]$ is a red solution of $L(m,a).$

 \vspace{.15in}
 
 \noindent\textbf{Case 3:}  $1\in R, \ 2\in B, \ 3\in R$
 
 \vspace{.15in}
 
 In this case  we will use numbers no larger than $3v.$  As in Case 2, all these numbers are in $[C(m,a)].$
 
 \vspace{.1in}
 
 \noindent\emph{Subcase 1:} $m$ is odd
 
 \vspace{.1in}
 
 First suppose $v$ is even. Then $c$ must be odd, so $c> 0.$ It follows that for any $k\in \{v+2,v+4,\ldots, 3v\},$ $ak$ is an even number such that $m-1\leq ak\leq 3(m-1),$ and therefore we can achieve the value $ak$ by assigning each variable on the left side of $L(m,a)$ a value of 1 or 3.  So $v+2,v+4,\ldots, 3v$ are all in $B$. To obtain a blue solution of $L(m,a),$ we start by assigning the value $3v$ to $x_m$ and the value 2 to each of the remaining variables. To achieve a solution, we must then increase the total value on the left side of $L(m,a)$ by $av+2-2c,$ which is easily seen to be at least $v$.  So we write $$av+2-2c=qv+r,$$ where $1\leq q\leq a, \ 0\leq r<v,$ and $r$ is even.  If we  increase the values of $x_1,\ldots,x_q$ to $v+2,$ we can then increase the value of $x_1$ to $v+2+r$ and obtain a blue solution of $L(m,a).$ Note that $v+2+r$ is even and at most $3v.$
 
 Now suppose $v$ is odd.  Then $v+1,v+3,\dots,3v-1$ are all even, and as in the preceding paragraph we see that they are all in $B$. To obtain a blue solution of $L(m,a),$ we start by assigning the value $3v-1$ to $x_m$ and the value 2 to each of the other variables. We must then increase the total value of the left side by $av+2-2c-a.$ Since $v$ is odd, $v\geq 3,$ and using this it is easy to show that $av+2-2c-a\geq v+1.$  So we write $$av+2-2c-a= q(v+1)+r,$$ with $1\leq q<  a,$ $r$ even, and $0\leq r\leq v-1$ since $v$ is odd.  If we increase the value of $x_1$ to $v+3+r$ and the values of $x_2,\ldots, x_q$ to $v+3,$ we obtain a blue solution of $L(m,a),$ since $v+3+r\leq 2v+2\leq 3v-1$ because $v\geq 3. $ (We could have done this argument by increasing values  to $v+1$ instead of $v+3,$ but doing it as we have will be useful in dealing with the next subcase.)
 
 \vspace{.1in}
 
 \noindent\emph{Subcase 2:} $m$ is even and $a$ is even
 
 \vspace{.1in}
 
 If $m$ and $a$ are even and $z$ is an even integer such that $m-2\leq (a-1)z\leq 3(m-2),$ then we can obtain a solution of $L(m,a)$ by assigning the value $z$ to $x_m$ and $x_{m-1}$ and assigning a value of 1 or 3 to each of the remaining variables. It is straightforward to verify that $m-2\leq (a-1)z\leq 3(m-2)$ whenever $v+2\leq z\leq 3v.$
 So if $v$ is even then $v+2,v+4,\ldots,3v$ are all in $B$, and if $v$ is odd then $v+3,v+5,\ldots,3v-1$ are all in $B.$
 
 We can now obtain a blue solution of $L(m,a)$ by repeating the arguments given for Subcase 1, because we didn't use the value $v+1$ in the argument given there when $v$ was odd.
 
 \vspace{.1in}
 
 \noindent\emph{Subcase 3:} $m$ is even and $a$ is odd.
 
 \vspace{.1in}
 
 In this subcase, $a-1$ is even, so we can now do the argument of the first paragraph of Subcase 2 without the restriction that $z$ be even, and conclude that $\{v+2,v+3,\ldots, 3v\}\subseteq B.$  We can then obtain a blue solution of $L(m,a)$ by using the argument given for even $v$ in Subcase 1, regardless of the parity of $v.$  In the present situation we will not know that the remainder $r$ is even, but that doesn't matter now.  $\Box$
 
 \vspace{.25in}
 
 \noindent\textbf{6. The proofs of Theorems 4 and 5}
 
 \vspace{.25in}
 
\noindent\emph{The proof of Theorem 4.}  Suppose that $\frac{2a}{3}+1\leq m\leq a.$ As in the proof of Theorem 2, we have $R_2(m,a)\geq 3,$ since $m\neq a+1.$

  It is shown in Theorem 2 of $[\textbf{6}] $ that $R_2(3,3)=9,$ so we can assume that $a\geq 4.$ If we take a bad 2-coloring of $[3]$ with $1\in R,$ then by statement (1) of Lemma 7 and the fact that $R_2(m,a)\neq 1,$ we have $2\in B.$ So by statement (3) of Lemma 7 we must have $3\in R.$  If $a\equiv m-1$ (mod 2), then by statement (2) of Lemma 7 we have a red solution of $L(m,a),$ so $R_2(m,a)=3.$

 If $a\not \equiv m-1$ (mod 2) then the coloring $R=\{1,3\}, B=\{2\}$ is bad, so $R_2(m,a)\geq 4$. To prove equality, suppose for a contradiction that we have a bad 2-coloring of $[4]$, with $1\in R.$ Then, as above, we have $2\in B$ and $4\in R.$  We again have $3\in R$ by statement (3) of Lemma 7. To obtain a red solution of $L(m,a)$ we assign the value 1 to all the variables, and show that we can increase the total value of the left side by $a-(m-1)$  by increasing some of the 1's on the left  side by 2 or 3 each.  This is possible if $2\leq a-(m-1)\leq 3(m-1).$  The second inequality holds since $m\geq \frac{a}{4}+1.$ The first inequality holds if $m\leq a-1.$ So we have a red solution unless $m=a.$ But if $m=a$ then, since $a\geq 4,$ $[a-4\rightarrow 3;\ 3\rightarrow 4;\ 1\rightarrow 3]$ is a red solution.  $\Box$
 
 \vspace{.15in}
 
 \noindent\textbf{Lemma 9.}  $R_2(4,5)=9.$
 
 \vspace{.15in}
 
 \noindent\emph{Proof.}  We first determine the unique bad 2-coloring of $[8]$ that has $1\in R.$  As in the proof of Theorem 4 we must have $2\in B$ and $4\in R,$ and it then follows from statement (2) of Lemma 7 that $3\in B.$ The solutions $[3\rightarrow 5;\ 1\rightarrow 3]$ and $[2\rightarrow 6;\ 2\rightarrow 3]$ then yield $5,6\in R,$ and the solutions $[2\rightarrow 7;\ 1\rightarrow 6;\ 1\rightarrow 4]$ and $[2\rightarrow 8;\ 2\rightarrow 4]$ show that $7,8\in B.$ With the coloring $R=\{1,4,5,6\}, B=\{2,3,7,8\},$ the left side of $L(m,a)$ would have total value at most 18 in any red solution, so $x_4$ would have to be assigned the value 1. But the left side couldn't have total value 5, so there is no red solution. In a blue solution, the left side of $L(m,a)$ would have total value at most 24, so $x_4$ would have to be 2 or 3. But the left side couldn't have total value 10 or 15, so there is no blue solution.
 
 Suppose now that we have a bad 2-coloring of $[9]$ with $1\in R.$ By what we have just shown, we must have $3\in B$ and $4,5,6\in R.$ Then the solution $[1\rightarrow 9;\ 3\rightarrow 3]$ shows that $9\in R,$ so $[1\rightarrow 5;\ 1\rightarrow 6;\ 1\rightarrow 9;\ 1\rightarrow 4]$ is a red solution.  We conclude that $R_2(4,5)=9.$ $\Box$
 
 \vspace{.15in}
 
 \noindent\emph{The proof of Theorem 5.} Suppose that $\frac{a}{2}+1\leq m <\frac{2a}{3}+1$ (so $a\geq 4$). Then the coloring $R=\{1\}, B=\{2,3\}$ of $[3]$ is bad by statement (3) of Lemma 7, so $R_2(m,a)\geq 4.$ For any bad 2-coloring of $[4]$ with $1\in R,$ we have $2\in B$ and $4\in R$ as above, so if $a\equiv m-1$ (mod 3) then by statement (4) of Lemma 7 we have a red solution of $L(m,a).$  So if $a\equiv m-1$ (mod 3) then $R_2(m,a)=4.$ 
 
Now suppose that $a\not \equiv m-1$ (mod 3).  Then, by statements (3) and (4) of Lemma 7, the coloring $R=\{1,4\}, B=\{2,3\}$ of $[4]$ is bad, so $R_2(m,a)\geq 5.$  Suppose we have a bad 2-coloring of $[5]$, with $1\in R.$  Then as above we have $2\in B$ and $4\in R.$  If $3\in R$ then, as in the second half of the last paragraph of the proof of Theorem 4, we have a red solution of $L(m,a)$ using values in $\{1,3,4\}.$ 
(The requirement $m\leq a-1$ at the end of the argument is satisfied since $m< \frac{2a}{3}+1$ and $m$ is an integer.)  So $3\in B.$
 
 We claim that $5\in R.$ To see this we obtain a solution of $L(m,a)$ using values in $\{2,3,5\},$ and note that any such solution must involve the value 5 by statement (3) of Lemma 7. To obtain our solution, we start by assigning the value 2 to all the variables.  We must then increase the total value of the left side by $2a-2(m-1)$ by increasing the values of some of the variables on the left side by 1 or 3 each.  As in the third-to-last paragraph of the proof of Theorem 2, to show that this is possible we need only verify that $0\leq 2a-2(m-1)\leq 3m-5. $ The first inequality holds since $m\leq a+1,$ and the second states that $2a\leq 5m-7,$ which is true since $m\geq \frac{a}{2}+1$ and $a\geq 4.$

 We now try to obtain a red solution of $L(m,a)$ by using values in $\{1,4,5\}.$ If we start by assigning the value 1 to all the variables, then to achieve a solution we must write $a-(m-1)$ as a sum of 3's and 4's, using at most $m-1$ terms.  This is possible if
 $ a-(m-1)\leq 4(m-1)$ (which clearly holds) and $a-(m-1)$ is either 3 or 4 or at least 6.
 Since $m-1< \frac{2a}{3}$ we have $a-(m-1)> \frac{a}{3},$ so $a-(m-1)\geq 2.$ So we have a red solution, and thus $R_2(m,a)=5,$ unless $a-(m-1)=2$ or 5.
 
 If $a-(m-1)=2$ then $m=a-1,$ so $\frac{a}{2}+1\leq a-1< \frac{2a}{3}+1$ and therefore $a=4$ or 5. For $a=5$ we have $R_2(4,5)=9$ by Lemma 9, and for $a=4$ it is shown in $[3]$ that $R_2(3,4)=10.$
 
 If $a-(m-1)=5$ then $m=a-4,$ so $\frac{a}{2}+1\leq a-4< \frac{2a}{3}+1$ and therefore $10\leq a\leq 14.$ It is easy to verify that in this case the coloring $R=\{1,4,5\}, B=\{2,3\}$ of $[5]$ is bad, so $R_2(m,a)\geq 6.$  For any bad 2-coloring of $[6]$ with $1\in R$ we have $3\in B,$ as above, and the solution $[5\rightarrow 6;\ a-9\rightarrow 3]$ shows that $6\in R.$ But then $[1\rightarrow 6;\ a-5\rightarrow 1]$ is a red solution of $L(m,a).$  So $R_2(m,a)=6.$  $\Box$

 \vspace{.25in}
 
\centerline{ \textbf{References}}
\vspace{.1in}

\noindent 1.  A. Beutelspacher and W. Brestovansky, Generalized Schur numbers, \emph{Lecture Notes in Mathematics}, Springer-Verlag, Berlin, \textbf{969} (1982), 30-38

\vspace{.03in}

\noindent 2. S. Guo and Z-W. Sun, Determination of the 2-color Rado
number for $a_1x_1+ \cdots +a_mx_m=x_0$,  \emph{J. Combin. Theory Ser.A}, \textbf{115} (2008), 345-353.

\vspace{.03in}

\noindent 3.  H. Harborth and S. Maasberg,  All two-color Rado numbers for $a(x+y)=bz,$ \emph{Discrete Math.} \textbf{197-198} (1999), 397-407.

\vspace{.03in}

\noindent 4. B. Hopkins and D. Schaal, On Rado numbers for
$\displaystyle \Sigma_{i=1}^{m-1} a_ix_i=x_m$, \emph{Adv. Applied
Math.} \textbf{35} (2005), 433-441.
\vspace{.03in}

\noindent 5. R. Rado, Studien zur Kombinatorik, \emph{Mathematische
Zeitschrift} \textbf{36} (1933), 424-448.
\vspace{.03in}

\noindent 6. D. Saracino, The 2-color Rado Number of $x_1+x_2+\cdots +x_{m-1}=ax_m,$  \emph{Ars Combinatoria} \textbf{113} (2014), 81-95.
\vspace{.03in}

\noindent 7. D. Schaal and D. Vestal, Rado numbers for $x_1+x_2+\cdots +x_{m-1}=2x_m$, \emph{Congressus Numerantium}  \textbf{191} (2008), 105-116.

\end{document}